\documentclass[a4paper,11pt]{article}

\usepackage{geometry}
\usepackage{authblk}
\usepackage{color}
\usepackage{xcolor}
\usepackage{graphicx}
\usepackage{hyperref}

\usepackage{pgf-pie}

\usepackage{bchart}
\usepackage{subfig}
\usepackage{setspace}
\usepackage{multirow}
\usepackage[square,numbers]{natbib}

\definecolor{coloraux1}{RGB}{211, 0, 119}
\definecolor{coloraux2}{RGB}{0,112,211}
\definecolor{coloraux3}{RGB}{104,223,234}
\definecolor{coloraux4}{RGB}{57,213, 137} 
\colorlet{color1}{coloraux1!95}
\colorlet{color2}{coloraux2!85}
\colorlet{color3}{coloraux3}
\colorlet{color4}{coloraux4!}

\usepackage{enumitem}
\usepackage[font={small}]{caption}

\begin{document}

\title{A Simple but Efficient Concept of Blended Teaching of Mathematics for Engineering Students during the COVID-19~Pandemic}
\author[1]{Saray Busto}
\author[1]{Michael Dumbser}
\author[1,2]{Elena Gaburro}

\affil[1]{Department of Civil, Environmental and Mechanical Engineering, University of Trento, Via Mesiano 77, 38123 Trento, Italy}
\affil[2]{INRIA, Univ. Bordeaux, CNRS, Bordeaux INP, IMB, UMR 5251, 200 Avenue de la Vieille Tour, 33405 Talence cedex, France}

\date{saray.busto@unitn.it, michael.dumbser@unitn.it, elena.gaburro@unitn.it}

	\maketitle
	
	\vspace*{-0.8cm}
	\begin{abstract}
		In this article we present {a case study} concerning a simple but efficient technical and logistic concept for the realization of blended teaching of mathematics and its applications in theoretical mechanics that was conceived, tested and implemented  at the Department of Civil, Environmental and Mechanical Engineering (DICAM) of the University of Trento, Italy, during the COVID-19 pandemic. The concept foresees \textit{{traditional}} blackboard lectures with a reduced number of students physically present in the lecture hall, while the same lectures are simultaneously made available to the remaining students, who cannot be present, via high-quality low-bandwidth \textit{{online streaming}}. {The case study presented in this paper was implemented in a single University Department and was carried out with a total of $n=1011$ students and $n=68$ professors participating in the study. }  %MDPI: Is italic necessary? ANSWER: we think so, it is used  throughout the paper to highlight the main points
	Based on our \textit{{first key assumption}} that traditional blackboard lectures, including the gestures and the facial expressions of the professor, are even nowadays still a very efficient and highly appreciated means of teaching mathematics at the university, this paper deliberately does \textit{{not}} want to propose a novel pedagogical concept of how to teach mathematics at the undergraduate level, but rather presents a \textit{{technical concept}} of how to \textit{{{preserve}}} the quality of traditional blackboard lectures even during the \mbox{COVID-19} pandemic and how to make them \textit{{available}} to the students at home via online streaming with adequate audio and video quality even at low internet bandwidth. The \textit{{second key assumption}} of this paper is that the teaching of mathematics is a dynamic \textit{{creative process}} that \textit{{requires}} the \textit{{physical presence}} of students in the lecture hall as audience so that the professor can instantaneously fine-tune the evolution of the lecture according to his/her perception of the level of attention and the facial expressions of the students. The \textit{{third key assumption}} of this paper is that students need to have the possibility to interact with each other personally, especially in the first years at the university.   
	We report on the necessary hardware, software and logistics, as well as on the perception of the proposed blended lectures by undergraduate students from civil and environmental engineering at the University of Trento, Italy, compared to traditional lectures and also compared to the pure online lectures that were needed as emergency measure at the beginning of the COVID-19 pandemic. {The evaluation of the concept was carried out with the aid of quantitative internet bandwidth measurements, direct comparison of transmitted video signals and a careful analysis of ex ante and ex post online questionnaires sent to students and professors.}
	\end{abstract}

% keywords 
	
	\section{Introduction and Context of This~Work}

The context of this work is the global COVID-19 pandemic, during~which many governments worldwide have imposed severe restrictions of all kinds of activities, including a complete shutdown of traditional lectures at schools and universities, which affected up to 83\% of the total enrolled learners in the world, i.e.,~almost 1.5 billion people~\cite{unesco2020education}. At~the University of Trento, Italy, regular lectures completely ceased on 5 March  2020 by Decree of the Rector, just after the beginning of the second semester (end of February 2020 to mid June 2020) of the academic year 2019/2020. As~a consequence, all teaching activities were entirely shifted to \textit{{online teaching}}, which was a common choice of many universities worldwide~\cite{chang2020countries}.   %MDPI: Is italic necessary? ANSWER: Idem (yes)

In view of the favorable pandemic development in Italy during the summer months of 2020, at~the University of Trento teaching activities with students present in the lecture hall were permitted again for the first semester of the academic year 2020/2021 (mid September 2020 to mid December 2020). However, due to national and regional COVID-19 restrictions, the~number of students allowed to enter each lecture hall was \textit{{reduced}} to 50\% of its nominal capacity; see the Decree of the Italian Prime Minister of  7 August~\cite{DPCM7Agosto}.	 
This restriction required the introduction of a suitable concept of \textit{blended teaching}, where some of the students are \textit{physically present} in the lecture hall and the others are enabled to follow the lectures \textit{online at home} via streaming over the internet. The~blended teaching concept proposed in this paper has a particular focus on those lectures that in pre-COVID times were traditionally held at the blackboard and that were usually highly appreciated by the students in the~past. 

The ideas outlined in this paper are therefore based on a \textit{first key assumption} that traditional blackboard lectures or chalk talks, which include the gestures and the facial expressions of the professor, are even nowadays still a very efficient means of teaching mathematics at the university level~\cite{FA12,AMea14,MWK17,WCWN19}. The~\textit{second key assumption} of this paper is that the teaching of mathematics is a dynamic and \textit{creative process} that \textit{requires} the \textit{physical presence} of students in the lecture hall as audience so that the professor can instantaneously fine-tune the evolution of the lecture according to his/her perception of the level of attention and the facial expressions of the students, in~addition to the questions asked and the comments made by the students. The~situation is similar to the one of an actor performing live in a theater in front of a big audience~\cite{AF11}. The~authors are convinced that this particular relationship between the professor and the audience not only improves the performance and the quality 
of the teaching for the audience present in the lecture hall, but~also the quality for those who follow the lectures online. In~the opinion of the authors, a~similar emotional experience is completely impossible when a professor teaches a pure online lecture merely in front of a computer screen, without~any possibility to directly see or perceive the emotions of his/her audience. 
The \textit{third key assumption} of this paper is that students need to have the possibility to \textit{interact} with each other personally~\cite{nerantzi2020use}, especially in the first years at the university when students should have an opportunity to get to know each other and to \textit{socialize}~\cite{Wen99}, although~this seems to be a contradiction to the social distancing measures adopted during the COVID-19 pandemic. 
Allowing the students to physically come to the lecture halls in order to live the lectures and exercises together as a common experience, despite the anti-COVID rules, seems to be still more attractive than following pure online lectures alone in front of a computer screen. For~a detailed analysis of digital readiness and the socio-emotional perceptions of higher education students see~\cite{Haendel}.      

As a consequence, our concept aims at enabling those students who can follow the lectures only online to have an experience that is as close to a real classroom experience as possible, permitting the students online to also communicate  directly with the professor and allowing them to have a high-quality view of the professor and the blackboard, combined with high audio quality. These are alltogether nontrivial requirements, in~particular due to pandemic \textit{internet bandwidth} restrictions. According to our three key assumptions, in this paper we therefore \textit{deliberately} do \textit{not} present a novel pedagogical concept of how to teach mathematics at the university, but~on the contrary, we aim at showing a simple but efficient \textit{technical concept} of how to \textit{preserve} the quality of traditional means of teaching mathematics at the university level in the context of the COVID-19 pandemic as far as~possible.

{
	% In this paragraph literature is inserted in a general context ... below wach of this literature voice are better described
	For what concerns the available literature on similar implementations, we can affirm that,
	as summarized in a recent review paper~\cite{ReviewPreprint} on face-to-face, blended and online teaching,  %MDPI: Is italic necessary? ANSWER: yes
	the definition of \textit{{blended learning}} is not unique.  
	Most of the times, the~term \textit{blended} refers generically to a combination of instructional methods, 
	which include (i) the availability of different online materials sometimes to 
	be studied independently by the learners before the lectures, in~the spirit of a flipped class strategy~\cite{CFW17,Cong2020},
	(ii) a mixture of online synchronous and asynchronous lectures~\cite{Galena2019-flipped,Ozadowicz}, 
	and, as~in our case, 
	(iii) a strategy to make possible the fruition of the same lecture from two different audiences, one in presence and the other one online.
	In all  cases, when a blended teaching concept is set up, a~central role is occupied by the use of new technologies 
	and its effectiveness and impact need to be investigated through dedicated~research.
	
	Innovative concepts of technology-based and \textit{blended teaching} had already been introduced in \textit{pre}-COVID times; 
	they were motivated by different logistic difficulties but all had similar objectives:
	implement, study and ameliorate new teaching strategies able to allow people to profitably attend lectures and to actively participate in them despite living in different places and with different time constraints.
	%2009
	A first remarkable example was given in~\cite{Kwan2009}, where the authors  analyzed and improved the teaching organization 
	in order to meet the need of \textit{working adult learners} that want to improve their knowledge and obtain new certifications but have limited free time. 
	Interesting findings of this research study were the preferences for synchronous and guided lectures and the request for a major interaction with colleagues and professors.
	%2016
	A second example is given by a \textit{joint} Spanish \textit{master} program in Industrial Mathematics between five different universities,
	%the Universities of Santiago de Compostela, A Coru\~na, Vigo, Carlos III of Madrid and Technical University of Madrid,
	which combines the best offer of lectures on applied mathematics of the involved institutions in order to provide a unique opportunity for talented students,  
	for which it would be, however, unfeasible to follow all these courses in presence in geographically distant campuses; for more details on the specific organization, see~\cite{M2I}.
	%2019
	A different blended approach is presented in~\cite{Galena2019-flipped},  
	where a new technological environment has been set up in order to allow various types of flipped classes, 
	i.e., to present and distribute introductory material (videos, quiz, documents) online 
	before a synchronous lecture takes place (in presence or online). 
	This research considers a \textit{pan-European network} of $18$ universities 
	and reports on the positive effects of the introduced strategy in activating students, 
	ameliorating their outcomes, increasomg flexibility and producing reusable material. 
	Moreover, in~\cite{gonzalez2019co}, the authors studied the effectiveness of a \textit{3D virtual campus}, 
	finding that postgraduate students really appreciate it for the increased flexibility, 
	while undergraduate students perceived it just as a complementary material. 
	For a similar investigation on blended approaches in the \textit{United States} 
	we refer to~\cite{McKnight}, where the authors analyzed six different strategies employed in US universities. %2016
	
	However, only with the \textit{advent} of the COVID-19 pandemic and the introduction of related restrictions 
	on movement and classrooms filling, blended teaching became central in many universities worldwide and, correspondingly, 
	a rich literature was developed both on the technical implementation of the concept and on the analysis of impact and satisfaction~rate. 
	
	During the first half of 2020 almost all university lectures were switched to a \textit{totally online} structure.
	In certain cases the initial reorganization was of emergency type and for example in~\cite{Naidoo2020}, %prima versione marzo 2020
	the authors focused on the necessity of increasing the professors' preparation on new technologies.
	Additionally, the authors of~\cite{Fan2021} emphasized  the role of teachers in this new and unexpected context.
	However, next to extemporary solutions, there are also many examples of \textit{well-organized} online courses, such as in~\cite{Ma2021}, where the authors claim the maintaining of excellent students' outcomes thanks to the exploitation of the Internet+ platform.
	Equally, in~the case-based paper~\cite{Fernandez2020}, the~authors report on the success of the online switch
	of a cybersecurity course that, already before the COVID-19 pandemic, was exploiting a blended approach  
	mixing in-presence classes and practical online computer laboratories. For an analysis of  online teaching related to COVID-19 at the University of C\'ordoba, Spain, see \cite{GHACM20}. 
	
	Together with the enthusiasm for the blended but completely online teaching solutions, 
	\textit{questions} on the social and cognitive implications of this approach started to arise.
	For example, in~\cite{2020_CFCovidTeaching} the authors present a review of articles written between January and April 2020
	analyzing the importance of a bidirectional communication and of keeping a certain degree of interaction between the students.
	Similarly, in~\cite{Gamage} we can find a description of the transformations necessary to teach  
	a large variety of STEM subjects online but also a careful discussion on the missing aspects of pure online teaching, 
	such as the reduced development of group cohesion and learners' communities, 
	as well as the increased inequalities due to different levels of internet access and private infrastructures. 
	Also in~\cite{Ozadowicz}, the~author describes the effective implementation of an automata course 
	that mixed synchronous and asynchronous online and offline methodologies, and~that, 
	while preserving remarkable students' outputs and increasing their effectiveness in autonomously searching and acquiring knowledge,
	lament the irreparable absence of a direct contact of engineering students with laboratory instruments. 
	On the contrary, according to~\cite{MDPIDecember2020Smith}, despite the reasonable alarm regarding  
	mental health symptoms in youth such as anxiety and depression (see~\cite{racine2020child}), 
	the outcomes of a sample of Canadian students was stationary during the pandemic and
	in particular extroverted students, in~the absence of other social stimulus, were more motivated to~study. 
	
	As part of the major criticisms and punctual analysis of pure online strategies, \textit{{drawbacks}}  %MDPI: Is italic necessary? ANSWER: yes
	have been highlighted in papers appearing at the end of 2020. 
	In particular the authors of~\cite{ReviewPreprint}, reviewing more than $70$ papers on different combinations of instructional methods,
	underline, among~many others, the~fact that lack of communication and social and physical contact provokes 
	diminishing oral capabilities and teamwork skills.
	Moreover, in~the conceptual work in~\cite{2020_PP_BlendedHigherEduactionCOVID}, 
	the authors highlight the necessity of social learning, peer-to-peer interaction and communication with professors,
	even in a sophisticated blended framework with a careful use of synchronous, asynchronous and flipped strategies. 
	Furthermore, both the last two cited works remark that pure online teaching considerably reduces 
	the possibilities for teachers to provide the needed inputs in time, support students and fine-tune the lecture depending on their reactions.
	Finally, in~\cite{2020_SR_COVIDGerman}, the authors analyze a case study in a German university focusing 
	on the problems connected with a lack of a dedicated technical equipment and the importance of socialization, 
	proposing some interesting solutions based on a hybrid campus~strategy.
	
	It is in view of these valuable studies, and~with the aim of mitigating the negative effects of pure online teaching while still respecting the imposed sanitary restrictions on in-presence lectures, 
	that we have designed, implemented and analyzed the blended teaching concept presented in this paper. 
	In particular, the~following objectives, questions and methodology have been considered:}

%MDPI: Is bold necessary? ANSWER: yes
{\textbf{{Purpose of this work.}} The purpose of this paper is to present the results of a case study about the design and implementation of a simple but efficient technical concept for blended teaching that aims at improving the quality of mathematics teaching at the university level during the adverse conditions imposed by the COVID-19 pandemic. In~particular, we want to show a simple technical solution to the rather challenging  problem of how to transmit complex mathematical blackboard content as well as audio signals at high quality via online streaming with limited internet bandwidth during pandemic conditions {{({note that during the COVID-19 pandemic the transmission of standard video signals via ZOOM and other online video conferencing software was limited to 360 p. Furthermore, overall internet bandwidth was essentially saturated due to compulsory telework in most public and private institutions})}}. As~such, the~concept presented in this paper may be useful also for other departments and universities with similar technical problems.% during the COVID-19 pandemic.    
	%MDPI: We changed footnote format throughout paper, please confirm. ANSWER: ok
	
	\textbf{Research questions.} The research questions of this paper are the following: Is it technically possible to realize high-quality video streaming of complex mathematical blackboard content combined with high-quality audio streaming during COVID-19 pandemic conditions with internet bandwidth limitations? What is the level of appreciation of the proposed technical concept for blended teaching from a student perspective? What is the level of appreciation and what are the encountered difficulties of the proposed blended teaching concept from a professor's point of view?

	\textbf{Research methodology and data analysis.} Concerning the adopted research methodology, a~carefully-designed experimental case study was carried out at the scale of an entire engineering department at the University of Trento, Italy, with~a total of $n=1011$ students and $n=68$ professors participating in the~experiment.

	In order to prepare the case study, the~necessary technical hardware and software were first of all very \textit{carefully selected} and \textit{thoroughly tested on-site} concerning technical performance and reliability at a small scale in a single lecture hall before scaling up the purchase of equipment to the entire department level with 18 lecture halls. In~particular, the~\textit{technical feasibility} of the concept was \textit{empirically verified} beforehand with the aid of \textit{quantitative internet bandwidth measurements} and with the aid of a detailed analysis of the obtainable video quality via a \textit{direct comparison of screenshots} taken under different transmission conditions. Indeed, in~this paper we show a direct comparison of several samples of transmitted mathematical blackboard content that can be achieved via standard video streaming and via the concept proposed in this paper. To~have a clear comparison with a reference solution, we also provide high-resolution photographs of the blackboard taken directly in the lecture~hall. 
	
	Furthermore, in~order to get reliable a priori information about the \textit{desiderata} of the students concerning blended teaching, an~ex ante online questionnaire was sent to all students. The~number of received answers for this first ex ante questionnaire (AS1) was $n=445$, which was a sufficiently large number to guide the design and conception phase of this case study and to get a clear idea about the students' preferences and expectations about blended teaching. Another ex ante online questionnaire (AP1) was sent to the professors, in~order to get a \textit{binding} answer concerning their preferred way of teaching, namely either \textit{fully online}, or~\textit{blended teaching}, with~lectures held \textit{live on-site} in the lecture hall and the \textit{simultaneous} transmission of these lectures over the internet \textit{via live online streaming} to the students at home, according to the technical concept presented in this paper. The~number of answers received from the ex ante questionnaire sent to the professors (AP1) was $n=68$. Right before the beginning of the semester, in~order to organize and implement the strict access rules and social distancing measures within the lecture halls, it was necessary to quantify the exact number of students who wanted to attend the blended lectures \textit{on-site} and of those who instead wanted to follow the blended lectures always \textit{online}. For~this purpose, a~second ex ante questionnaire (AS2) was sent to all students, receiving a total number of $n=1011$ answers.

	In order to assess the level of appreciation of the proposed technical concept by the students and by the professors, a~detailed ex post analysis based on two further online questionnaires was carried out. In~particular, with~the ex post questionnaire sent to the students (PS1), which received a total number of $n=509$ answers, the~two main questions concerning the perceived audio and video quality of the technical concept proposed in this paper were investigated more thoroughly, but~also a fairly generic overall level of satisfaction with the proposed blended teaching methodology was measured. The~ex post questionnaire sent to  the professors from the area of mathematics and its application to theoretical mechanics (PP1) received a total number of $n=6$ answers. 
	At this point we would like to stress again that the main focus of this paper is indeed on the blended teaching of \textit{mathematics and its applications}, due to the particularly complex blackboard content that needs to be transmitted. This explains the small sample size in the answers to the ex post questionnaire PP1. A~short summary of the different questionnaires used for the design and assessment of the present case study can be found in Table~\ref{tab.quest}. 
}

\begin{table}[] 		
	\caption{{   
			Summary of the different \textit{ex ante} (A) and \textit{ex post} (P) questionnaires sent to the students (S) and to the professors (P) and which were employed for the design, conception and assessment of the present case study.  } 
	}
	{\small
		\begin{tabular}{llll}
			\hline 
			questionnaire           & main purpose of the  questionnaire                                & type               & number $n$ of    \\ 
			number                  &                                                                   &                    & answers received \\
			\hline 
			\multirow{2}{*}{AS1}    & obtaining students' preferences and expectations                  & \multirow{2}{*}{\textit{ex ante}} & \multirow{2}{*}{\textbf{445} students}  \\ 
			& concerning blended teaching to design the case study              &                                   &                                \\[3pt] 
			\multirow{2}{*}{AP1}    & binding choice of the professors whether to teach                 & \multirow{2}{*}{\textit{ex ante}} & \multirow{2}{*}{\textbf{68} professors}   \\ 
			& onsite (blended), or purely online                                &                                   &                                \\[3pt] 
			\multirow{2}{*}{AS2}    & binding choice of the students to attend onsite or online         & \multirow{2}{*}{\textit{ex ante}} & \multirow{2}{*}{\textbf{1011} students} \\ 
			& needed for the strict access rules to the University buildings    &                                   &                                \\[3pt] 
			\multirow{2}{*}{PS1}    & evaluation of the students' appreciation of the proposed          & \multirow{2}{*}{\textit{ex post}} & \multirow{2}{*}{\textbf{509} students}  \\ 
			& blended teaching concept implemented in this case study           &                                   &                                \\[3pt] 
			\multirow{2}{*}{PP1}    & evaluation of the professors' opinion concerning the proposed     & \multirow{2}{*}{\textit{ex post}} & \multirow{2}{*}{\textbf{6} professors}    \\ 
			& concept for the area of mathematics and theoretical mechanics     &                                   &                                \\ 
			\hline 	     
		\end{tabular} 
	}
	\label{tab.quest} 
\end{table}

The rest of this paper is organized as follows: in Section~\ref{sec.concept} we present and motivate the concept of blended teaching developed and implemented at DICAM during the first semester of the academic year 2020/2021. The~motivation of our blended teaching concept is also  in  light of the \textit{pure online} teaching that was adopted as a compulsory \textit{emergency measure} during almost the entire second semester of the academic year 2019/2020. In~Section~\ref{sec.realization} we describe the details of the technological and logistic realization of the concept, in~particular  the hardware and software requirements, the~special training of the professors before the start of the semester and the organization of the technical support for blended teaching that was made available during the semester. The~quality of the adopted concept was quantitatively evaluated by both professors and students via ex post online questionnaires. The~results of this evaluation are presented and discussed in Section~\ref{sec.evaluation},  
{which also contains a direct \textit{quantitative validation} of the low-bandwidth high-quality screen sharing approach proposed in this paper by showing a set of comparative screenshots of the new approach versus standard video 
	streaming, together with quantitative measurements concerning the required internet bandwidth.} Finally, in~Section~\ref{sec.conclusions}, we propose some concluding~remarks. 

\section{Blended Teaching Concept at DICAM during the COVID-19~Pandemic} 
\label{sec.concept}

In order to illustrate and to motivate the blended teaching concept elaborated and implemented by DICAM during the first semester of the academic year 2020/2021 (\mbox{Sections~\ref{sec.MainPillars}} and \ref{sec.CovidPillar}), we first briefly describe the typical structure of the courses in pre-COVID times (Section~\ref{sec.PreCovid}) and also draw some conclusions from the purely online teaching phase that was mandatory in spring 2020 at the beginning of the pandemic (Section~\ref{sec.EmergencyPhase}).

\subsection{Typical Structure of the Courses in Mathematics Offered at~DICAM} 
\label{sec.PreCovid}

All courses at DICAM in the areas of mathematics and its applications, including theoretical fluid mechanics and solid mechanics, which are under consideration in this study, follow the traditional scheme of theoretical lectures combined with classroom exercises. In~addition, for~courses offered at the bachelor level, there are special group exercises organized in small work groups, led by MSc or PhD students (tutors). The~group exercises are based on the elaboration and discussion of exercise sheets in small work groups and make use of a flipped-class concept, where students first elaborate the exercises at home and then discuss the results with the tutors in the classroom~\cite{Mac15,CFW17}. Online surveys and structured interviews with  student representatives held in pre-COVID times revealed a clear preference of the students of DICAM for traditional blackboard lectures, documented in the \textit{Yearly report of the joint committee of professors and students of DICAM}, which is not available for public~view.    

\subsection{Experience of the Emergency Phase during the Second Semester 2019/2020}  
\label{sec.EmergencyPhase}

With the shutdown of regular lectures in March 2020, all teaching activities at DICAM suddenly needed to be held entirely online. In~addition to the combined online teaching platform Moodle plus 
Kaltura {(\url{https://moodle.org} and \url{https://corp.kaltura.com/})}, which was already available in pre-COVID times, for~the emergency online teaching the University of Trento also provided its professors with licenses for the software ZOOM {(\url{https://zoom.us})}. Online teaching could be held either in a \textit{synchronous manner} by direct online streaming of the lecture content via ZOOM, or~in an \textit{asynchronous way} by pre-registering the lectures and uploading them into the combined Kaltura--Moodle platform. The~personal experience of the second author of this paper  with  pure online teaching activities from March to June 2020 was overall rather negative, mainly due to the complete lack of an instantaneous feedback from the students (missing facial expressions, comments and questions) even during synchronous teaching, since all students kept their microphones and webcams systematically switched off in order to save internet bandwidth, but~also due to the lack of the blackboard as a traditional means of teaching mathematics at the university. 
%Hence, obtaining an instantaneous feedback of the students during the synchronous online lectures was overall rather difficult. Also the overall level of satisfaction of the students with the emergency teaching online was only average, as revealed by systematic online surveys of the Department. 
In a systematic survey {(AS1)} made by the department in June 2020, at~the end of the emergency online teaching during the second semester of 2019/2020, {ex ante to the case study presented in this paper} and to which {\emph{n}~=~445} students of DICAM responded, the~\textit{clear preference} of the students regarding online teaching was  \textit{synchronous online lectures} with the recording of the lecture \textit{made available} to the students via the Moodle--Kaltura platform afterwards. In~the same {ex ante} online survey {(AS1)} the students also expressed a \textit{clear preference} for \textit{traditional lectures at the blackboard}, confirming that the body language of the professor is important and that the possibility to see the professor renders the lecture more interesting. For~a summary of the main results of this survey, see Table~\ref{tab.as1.survey}, from~which it also becomes evident that the majority of  students did not consider the option to attend lectures only online in the case a traditional lecture would again be possible in the subsequent~semester.

 \begin{table}[] \small
	\caption{\textcolor{black}{\textit{Ex ante} survey (AS1) } of DICAM made among its students in June 2020 at the end of the emergency online teaching  of the second semester of 2019/2020. Total number of students participating in the survey: $n=445$. Not all questions were mandatory.}
	\label{tab.as1.survey} 
	\begin{tabular}{ll}  
		\hline 
		\textbf{Questions} & \textbf{Answers with absolute numbers} \\ 
		\hline 
		1. Which online teaching modality do you prefer? &  asynchronous (47) \\
		&  synchronous without registration (17) \\
		&  \textbf{synchronous with registration (257)} \\ 
		&  depends on the lecture (67) \\
		&  depends on the professor (42) \\
		&  no opinion (15) \\
		\hline 
		2. Seeing the body language and the facial  	 &  0 - no opinion (20) \\
		expressions of the professor makes the  		 &  1 - not at all (1) \\
		lecture more interesting? 					 &  2 - a bit (17) \\ 
		&  \textbf{3 - a lot (87)} \\ 
		\hline 
		3. The use of the blackboard as a traditional  	 &  0 - no opinion (29) \\
		means of teaching is useful?   		 			 &  1 - not at all (1) \\
		&  2 - a bit (11) \\ 
		&  \textbf{3 -  a lot (81)} \\ 
		\hline 
		4. Do you consider the option to attend the lectures     &  \textbf{0 - no (263)} \\
		of the first semester of 2020/2021 only online?   		 & 1 - yes (162) \\
		\hline 
	\end{tabular}
\end{table}

Based on the results presented in Table~\ref{tab.as1.survey} and in order to overcome the shortcomings of the purely online emergency teaching of the second semester 2019/2020, and taking also into account the updated national and regional regulations that permit regular teaching activities under certain conditions, the~Department of Civil, Environmental and Mechanical Engineering (DICAM)  adopted the  strategy for the first semester of the academic year 2020/2021 detailed in the next~section.    

\subsection{Main Pillars and General Objectives of the Blended Teaching~Concept}
\label{sec.MainPillars}

Given the clear preference expressed by the students of DICAM in favor of \textit{traditional blackboard lectures} actually held in the lecture hall and given also the reciprocal importance of the body language and the facial expressions of the professor seen by the students and, vice~versa, the~ones of the students seen by the professor, the~department has therefore decided to adopt a simple blended teaching concept that tries to create an environment that is \textit{as normal as possible}, both for~the professors and for the~students. 

The concept consists in the \textit{high-quality transmission} of traditional blackboard lectures that are actually held in the lecture halls in front of up to 50\% of the students who are allowed to be physically present and to make this experience \textit{available} as much as possible via \textit{low-bandwidth online streaming} to the remaining students who have to follow the lectures online at home due to  COVID-19 restrictions. The~two most important pillars on which our blended teaching concept relies are the \textit{high-quality view} of the blackboard and the professor, combined with \textit{high audio quality} to capture the speech of the professor, not only for the students present in the lecture hall, but~also for those online. Under~pandemic internet bandwidth restrictions, both pillars were anything  but trivial to realize, but~in the next section we will give technical details of how these objectives can be reached with standard hardware and software.   
Furthermore, \textit{bidirectional communication} is an integral part of the blended teaching concept discussed in this paper, thus allowing  the students who follow the lectures online to communicate \textit{directly}  with the professor, and~making their contributions, questions and comments also available to all of the students who are physically present in the lecture hall.   
As already stated in the introduction, the~blended teaching concept adopted by DICAM \textit{deliberately} does \textit{not} present a novel \textit{pedagogical concept} of how to teach mathematics, but~on the contrary, aims at demonstrating a simple but efficient \textit{technical concept} of how to \textit{preserve} the quality of traditional means of teaching mathematics in the context of the COVID-19 pandemic as much as possible. 
The three main objectives were: 
(i) To \textit{preserve} the high quality of traditional blackboard lectures that were usually appreciated by the students in the  pre-COVID era also during the COVID-19 pandemic and to make this form of lecture also \textit{accessible} to all students who were not able to attend the lectures in person, in~particular also to those with \textit{low internet bandwidth};   
(ii) To allow students to interact personally with each other and with the professors, complying with all of the restrictions due to the COVID-19 pandemic such as maintaining a minimum distance between each other and the obligation to wear a face mask during the entire period in which the students are present in the building;   
(iii) To allow professors to continue using their well-established traditional teaching concepts as much as possible, even under pandemic conditions.  
The same blended teaching concept was also adopted for the group exercises led by the tutors, allowing small work groups of students to attend the group exercises in person, while the others were following~online. 

The obvious \textit{shortcomings} of the blended teaching concept illustrated in this section are the rather stringent requirements on the computer hardware and software to be used in the lecture hall, combined with the necessary technological skills that needed to be acquired by the professors before and at the beginning of the~semester. 

To realize the concept detailed above, \textit{all} lecture halls of the department needed to be technically upgraded with an appropriate audio and video system, since no special equipment for blended teaching was present at DICAM in the pre-COVID era. The~department furthermore organized a series of special training sessions, in~which the professors were trained to use the necessary hardware and software in order to get ready for blended teaching before the beginning of the semester. In~order to allow the students to get a preview of how the blended teaching would be in the first semester of 2020/2021 and in order to test the equipment and to show professors the technical possibilities, the~students were invited to join some of these training sessions~online.     

\subsection{COVID-Specific Part of the~Concept} 
\label{sec.CovidPillar}

In order to establish a priori which students were allowed to enter the lecture hall in a given week and~which were not, the~department  carried out {another ex ante} online survey { among the students (AS2)} before the beginning of the semester, in~which the students could express their intention whether they wanted to attend lectures in person or~not. The~choice was \textit{free}, but~\textit{binding}. {In total, $n=1011$ students provided answers to (AS2).} As a result of this survey (see Table~\ref{tab.as2.survey}), in~which {$k=181$ and therefore only $k/n = 18\%$} of all students chose to follow lectures exclusively online, the~department divided those students who wanted to attend lectures personally into sub-groups that would be allowed to come to the lecture halls based on a \textit{weekly rotation principle} so that all students had the same number of weeks in which they could attend the lectures in person and in which they needed to follow the lectures online. A~maximum nominal capacity of 50\% of each lecture hall combined with the percentage of students who decided to attend all lectures only online allowed most of the students of DICAM to be  present in the lecture halls for most of the~weeks.    %unclear what is meant by "most of the weeks".

\begin{table}[]
	\caption{\textcolor{black}{\textit{Ex ante} survey (AS2) of DICAM made among its students in September 2020, right before the start of the lecture period. Needed for the top-down access authorizations to the University buildings and the organization of the groups of students allowed to attend blended lectures onsite according to a rotation principle. Total number of students participating in the survey: $n=1011$. } } 
	\label{tab.as2.survey} 
	\textcolor{black}{ \small
		\begin{tabular}{ll} 
			\hline 
			\textbf{Question} & \hspace{-2.5cm} \textbf{Answers with absolute numbers} \\ 
			\hline 
			In the next semester, do you want to attend blended lectures onsite?  &  \textbf{yes (830)} \\
			Note: if you choose \textit{no}, you will \textit{not} be authorized  &   \multirow{2}{*}{no (181)} \\
			to enter the University buildings during this semester. & \\
			\hline 
		\end{tabular}
	} 
\end{table}

An important choice of the adopted concept was to allow each professor to use his/her own laptop for the blended teaching, for~two reasons: 
(i) Using their own laptops, professors can use their preferred computer  environment and operating system with which they have the most experience and which they are most acquainted with;  
(ii) Using their own laptops avoids touching common computer keyboards and pointing devices and therefore solves the problem of disinfection of these devices between one lecture and the~other.       

To reduce the stress induced by the obligation of wearing face masks continuously while present in the university building and to reduce the pressure on the public transport system, the~timetable was re-organized in a morning block (8:30--13:30) and an afternoon block (14:30--19:30), so that some students were present only in the morning block and others only in the afternoon block. As~an anti-COVID measure the department allocated one and the same lecture hall only for at most two different groups of students, with~an interval of one hour left between the two groups that was needed for cleaning purposes. Furthermore, one medium-sized lecture hall was kept empty as a strategic reserve in case of COVID-related disinfection measures or in case of technical difficulties in another lecture hall. The~adopted timetable leaves a total amount of at most 25 h of real lectures in the lecture hall per group, which is not always sufficient to fit the entire timetable of each group of students. To~mitigate the effects of this decision, before~the start of the semester, all professors of the department {had to answer the ex ante questionnaire (AP1) (see \mbox{Table~\ref{tab.ap1.survey}}), in~which they} could choose in a free but binding manner whether to offer their lectures purely online, or~whether they wanted to adopt the blended concept outlined above. {The total number of answers received for questionnaire (AP1) was $n=68$ and the number of professors opting for blended teaching was $k=47$.} Therefore, most professors preferred the blended option for their courses {($k/n=69\%$)}, but~the remaining 31\% who chose the online teaching made it possible to fit all blended courses into the given 25 h limit. The~pure online lectures were mostly held in synchronous manner and could be followed by the students \textit{only at home} and \textit{not} from the university building. For~this reason, the~pure online lectures were scheduled only at a 2 h distance from blended lectures, in~order to allow students to move from the university building to their homes and~vice~versa. {In survey (AP1) the professors were also asked whether they agree with the proposed staggered timetable. $k=62$ and therefore an overwhelming majority of $k/n=91\%$ of the professors answered this question with \textit{yes}, thus allowing  the department to proceed with the organization of the timetable as outlined above. }

\begin{table}[!ht]
	\caption{\textcolor{black}{\textit{Ex ante} survey (AP1) of DICAM made among its professors in June 2020 for the organization and implementation of the blended teaching concept presented in this paper. Total number of professors participating in the survey: $n=68$. } } 
	\label{tab.ap1.survey} 	\vspace{-5pt}
	{\small 
		\begin{tabular}{ll} 
			\hline 
			\textbf{Questions} & \hspace{-1.1cm} \textbf{Answers with absolute numbers} \\ 
			\hline 
			1. How do you plan to teach your course next semester? &  \textbf{onsite blended (47)} \\
			\hspace{2.3mm} Note: if you choose \textit{exclusively online}, no lecture hall will be assigned &  \multirow{2}{*}{exclusively online (21)} \\
			\hspace{2.3mm} to your course for the entire semester & \\ 
			\hline 
			2. Do you agree with the proposed staggered timetable  & \textbf{yes (62)} \\
			\hspace{2.3mm} between blended onsite and purely online lectures?                                           &  no (6) \\
			\hline 
		\end{tabular}
	} 
\end{table}         

 The entrance and the exit of the university buildings were controlled electronically via a smartphone app that was specifically developed by the ICTS %define if appropriate.
services of the University of Trento and which required the scanning of a QR code upon entry and exit of the university buildings. The~groups of students who were allowed to attend lectures in a given week were authorized top--down by an electronic authorization system that was {based on the choice expressed by the students in the ex ante survey (AS2)} and was coupled with the cellphone app. Queues at the entrance were reduced by slightly shifting the start of the lectures at the MSc level and those at the BSc level by 15~minutes.

\section{Technological and Logistic Realization of the~Concept} 
\label{sec.realization} 

To realize the concept described above, the~personal laptop of each professor became the crucial hub for the acquisition and distribution of the audio and video streams to be sent online. As~a videoconferencing system we used ZOOM, already licensed by the University of Trento during the first phase of the pandemic and due to its very good performance concerning audio and video live streaming even in the case of low available internet bandwidth.    
In order to minimize technical problems a priori, two completely \textit{redundant} systems were installed in each lecture hall: a~videocamera system that allows to transmit traditional blackboard lectures, but~also \textit{classroom experiments}, as~well as a graphics tablet and digital pen system~\cite{WM13,Mac14}.     		
In this manner, each professor was also offered the free choice of which means he/she wanted to use for teaching, possibly also a combination of both (videocamera + tablet). The~videocamera system was also mandatory for those lectures that involve classroom experiments, which is very frequent in theoretical solid mechanics~\cite{Bigoni}, see,  
e.g., {(\url{https://bigoni.dicam.unitn.it/})}.  
The resulting main \textit{challenges} that needed to be overcome to implement this strategy were the~following:
\begin{itemize}
	\item The lecture halls at DICAM were not at all prepared for blended teaching and thus needed to be specifically equipped;  
	\item Achieving sufficiently high video resolution online that allows to read complex mathematical formulas clearly, despite the existing COVID-19 bandwidth and HD video feature restrictions;%, was the most difficult technical problem to overcome;   	 	  
	\item Achieving at the same time a sufficiently high audio quality both in the lecture hall and online, %was challenging
	in particular capturing the speech with high quality while the professor is turned towards the blackboard while writing and allowing also a clear bidirectional communication;    
	\item %Another major difficulty was due to 
	The fact that rather complex and relatively new technology needed to be handled by each professor during his/her blended lectures;  
	\item The highly heterogeneous computer hardware and operating systems resulting from the choice of allowing each professor to use his/her personal laptop, which also 		
	made the training and technical support particularly challenging and  led to a wide set of necessary adapters to connect to the standard hardware installed in each lecture hall. 
\end{itemize}

In order to convey an idea of the scale {of the present case study}, let us also note that the concept presented in this paper was applied to the \textit{entire} department {with 18 operational lecture halls}, covering more than 50 courses from the fields of mathematics, physics, chemistry, engineering and social sciences, {including 1011 students and 68 professors}. In~what follows, we will further comment on the above challenges, indicating the strategies employed to surmount~them.

\subsection{Description of the Employed~Equipment}

The implementation of the developed blended teaching concept required the lecture halls to be furnished with a completely new set of electronic devices. The~initial equipment of the classroom accounted for a blackboard, the~projector system and, in~large lecture halls, also a sound system so that the professor can be clearly heard on site. On~top of this we introduced:

\begin{itemize}
	\item Personal laptop. As~already mentioned, the~particular setup designed in this case study gives the professor the opportunity to employ his/her own laptop, aiming at minimizing the novelties resulting from the use of new technical devices. 
	Consequently, the~laptop of the professor became the main control device of the blended teaching~lectures.
	
	\item Camcorder. To~properly capture the blackboard and transmit physical experiments, a full high-definition camcorder was needed. It was further equipped with a directional microphone to record the speech of the professor in small classrooms. For~large blackboards, instead of the standard camcorder proposed here, it was necessary to employ a wide-angle lens camera (e.g., PTZ Minrray {UV540}) controlled by IR %define if appropriate.
	remote control to allow switching between different parts of the blackboard and different levels of~zoom.
	
	\item Camlink 4K. The~transmission of a high-quality video signal from the camcorder to the computer required  a 4K HDMI to USB converter. Let us note that although the quality broadcast by these devices goes up to 1080 p at 60 fps or 4K at 30 fps, to~reduce bandwidth consumption we suggest to limit the video capturing program to 1080 p at 25 fps. 
	%		It is an HDMI tethering device required for the transmission of a high quality stream video signal from a camcorder (that has an HDMI output) to the computer. It results in an HDMI to USB converter with broadcast quality up to 1080p at 60 fps or 4K at 30 fps.		
	% Para mi thetering se usa solo cuando la conexion es inalambrica, no fisica, nosotros no transmitimos la clase con estas resoluciones-fps sino mas bajas (gestionado por VLC/OBS/Quick): 1080p 30fps max por mucho que el dispositivo lo pase al ordenador con esta calidad, ¿podria dar lugar a equivocos si ponemos esto y no avisamos?

	\item Active USB 3.0 extension cable. Typically the camcorder and, consequently, the~camlink are far away from the laptop, and thus a USB 3.0 extension cable is needed. To~properly preserve the video signal and provide the needed power supply for the camlink, an~active cable was~employed.
	
	\item USB Hub 3.0. To~minimize the number of connectors to be plugged into the laptop, they were all gathered in a unique USB hub of four ports that transmits input and output signals to and from the computer. At a minimum, SuperSpeed USB (USB 3.0) cables, connectors and ports are needed to transfer all of the signals,  the use of USB 3.1 or 3.2 ports being preferred over USB 3.0~ports.
	
	\item Graphics tablet and digital pen system. The~alternative system to the blackboard was a rubber grip HD 22 inch graphics tablet that was not hand sensible (i.e., a pen display), so that the writing is as close as possible to traditional handwriting on a piece of paper, used together with the preferred interactive whiteboard software. To~be able to import and export PDF documents and save the lecture notes we  suggested the use of the OpenBoard {(\url{https://openboard.ch/index.en.html})} software.
	Communication of the graphics tablet with the computer was achieved via an input HDMI signal that duplicated the laptop screen and an output USB signal for the recording of the~handwriting. 
	
	\item HDMI splitter. The~laptop screen must be shared both with the students in the lecture hall using the projector system and with the graphics tablet. To~duplicate the signal an HDMI 4K splitter is employed. Since most projectors at DICAM are only equipped with a VGA input port, an HDMI to VGA converter is connected to the HDMI cable exiting form the HDMI~splitter.
	
	\item Airlink. To~facilitate bidirectional communication, the comments of online students were broadcast through the pre-existing sound system of the large lecture halls using an airlink connected to the laptop via bluetooth. For~small and medium classrooms (up to 100 students) and without a pre-installed sound system an alternative was the use of a speakerphone (Jabra Speak 810) that can be connected to the laptop using the audio jack or via bluetooth. The~latter device allows listening to the online comments in the lecture hall and also captures  the voice of the people present in the classroom to be transmitted online, hence guaranteeing peer interaction, a~key point of laboratories and seminars. Meanwhile, the~former device only allows the online students to be heard in the lecture~hall.
	
	\item Ear microphone. In~large lecture halls and when the professor is turned towards the blackboard, the~use of a directional microphone was not enough to properly capture the speech with high quality. To~solve this problem a personal bluetooth ear microphone was employed by each~professor.  
	
	\item Adapters. The~resulting hardware of each lecture room, as~explained below, requires the personal laptop of each professor to be connected thorough one HDMI and one USB/A port. Due to the great heterogeneity of modern computer hardware, a~large set of adapters was made available, allowing to also use  VGA, DisplayPort, mini DisplayPort, mini HDMI, micro HDMI and USB-C~ports.
	
\end{itemize}

A detailed list of the above components and their estimated cost can be found in \mbox{Table~\ref{tab.components.budget}}.

 \begin{table}	
	\caption{List of selected components employed in each lecture hall and corresponding cost estimate in EUR including 22\% VAT. The equipment was carefully selected and tested before installation. } 
	\label{tab.components.budget} 
	\begin{center} \small 
		\begin{tabular}{lr} 
			\hline 
			\textbf{Equipment} & \textbf{Cost estimate in EUR} \\ 
			\hline 	
			Laptop of the professor			&				0     \\ 
			Panasonic HC-V770 camcorder		&			    408   \\
			Elgato Camlink 4k HDMI to USB converter	&				130   \\ 
			Ugreen 10m active USB 3.0 extension cable 	& 	43    \\ 
			Anker USB hub 3.0				&				13    \\
			Wacom Cintiq 22 graphics tablet	&				953   \\ 
			4k HDMI splitter Ablewe			&				16    \\
			Primewire 2m HDMI cable		    &				10    \\ 				
			Ugreen HDMI to VGA adapter	    &				10	  \\
			BT DW 20BR Klark Teknik bluetooth airlink	&	89	  \\
			Jabra talk 25 bluetooth ear microphone{*} 	    &   28    \\		
			\hline 
			\textbf{Total amount:}   						  & \textbf{1700} \\ 
			\hline 
			\multicolumn{2}{l}{{*}one ear microphone was needed for each professor to comply with anti-COVID regulations					}  		  \\ 
		\end{tabular}
	\end{center} 
\end{table}

\subsection{Physical Setup of the~Equipment}

\begin{figure}[!t]
% 	\begin{center} 
\includegraphics[width=0.95\textwidth]{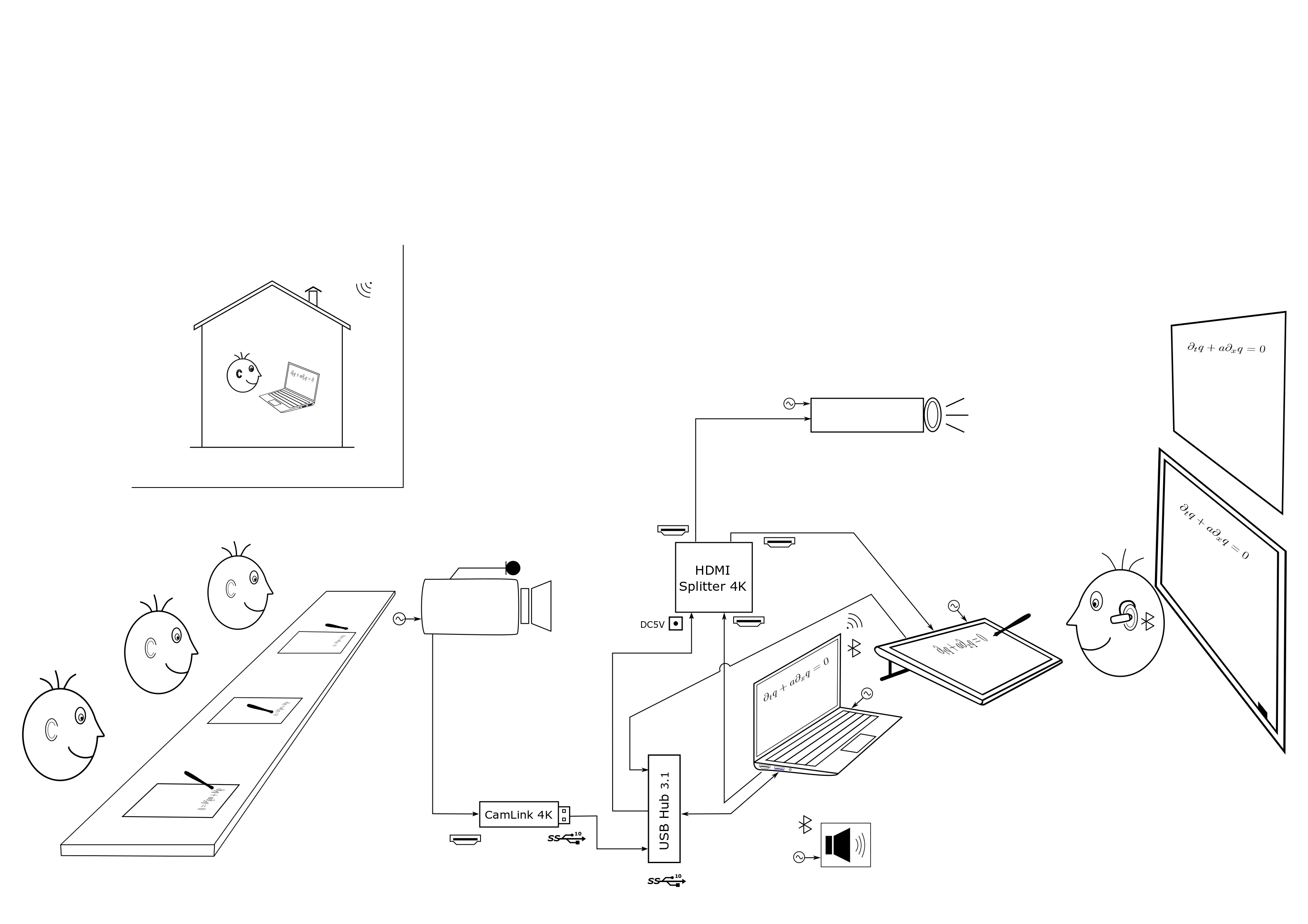} 		
% 	\end{center} 
\caption{Sketch of the technical realization of the blended teaching concept at~DICAM. } 		
\label{fig.scheme} 
\end{figure}

%	The needed technological devices to meet all the requirements for a high quality blended teaching are reported in Table~\ref{tab.components.budget}. 
A sketch of the connections among the different devices is shown in Figure~\ref{fig.scheme}. 	
On the one hand, to~allow high-quality capturing of the blackboard, the~full HD video camera is connected via HDMI cable to a camlink 4K, which allows the transition to a USB 3.1 signal transmitted using an active USB extension cable to a USB hub connected to the laptop of the professor. 	
The alternative and/or complementary device for the blackboard is the pen display, which transfers the data collected via a USB 3.1 cable connected to the USB Hub. Meanwhile, the~signal of the computer arrives at the tablet using an HDMI cable. To~allow the display of the laptop screen both in the graphics tablet and in the classroom projectors an active HDMI splitter 4K is employed. Consequently, a~single HDMI cable will connect the laptop to the HDMI splitter, from where two different HDMI cables transmit the same signal to the tablet and the projector; if necessary, an HDMI to VGA adapter is employed. Let us note that the HDMI splitter is also connected to the USB hub from which it receives the needed power supply. 	
The speech of the professor is picked up from a bluetooth ear microphone and collected in the computer. To~ease bidirectional communication the laptop is connected using a bluetooth airlink to the sound system of the lecture~hall.
% start a new page without indent 4.6cm
%\clearpage

\nointerlineskip

\vspace{-6pt}

\subsection{High-Quality Audio and Video Online~Streaming}

%MDPI: Is bold necessary? ANSWER: yes
\textbf{{High}-quality low-bandwidth video streaming via screen sharing.}  
In order to share visual information, the~lecture halls were equipped with full HD video cameras that capture the blackboard content and~graphics tablets.	
In-presence students can see the contents directly on the blackboard or on the projector screens. 
For online students the visual information was streamed via ZOOM.	
%	 In addition, the lecture halls were equipped graphics tablets and digital pens, as well as HDMI splitters to connect the laptop of the professor and the graphics tablet to the existing video projectors in the lecture hall.  
%	
However, ZOOM automatically reduces the video quality to $640\times360$~p due to COVID-19 internet bandwidth restrictions for all meetings with three or more people simultaneously connected. This resolution was \textit{far too low} to transmit the content of a blackboard with complex mathematical formulas, containing  small indices and superscripts, properly to the online students. The~only option to overcome this intrinsic limitation and to obtain a proper transmission of the video signal of the blackboard in permanent full HD resolution ($1920 \times 1080$~p) was to use the \textbf{share screen feature} of ZOOM. While this was rather obvious in the case of the use of a graphics tablet, it required the use of a third party video capturing software to capture and transmit the signal coming from the video cameras installed in the lecture hall.  
We found the following free video capturing softwares particularly useful: 
VLC and OBS ({\url{https://www.videolan.org/} and \url{https://obsproject.com/}}) on Windows and Linux, and Quicktime ({\url{https://www.apple.com}}) on Macintosh computers. The~reasonable \textit{compromise} in this choice lies in the \textit{reduced framerate} of the video, while~maintaining the full HD resolution at all times, which is \textit{crucial} to transmit the static content of the blackboard in optimal~quality.

\textbf{Audio streaming.} 
The bluetooth ear microphone is \textit{crucial} for blackboard lectures in order to suppress background noise and allow the professor to talk while writing on the blackboard. This ear microphone, which serves to transfer the audio in high quality to the students attending the lecture online, is worn \textit{in addition} to the standard classroom microphone that serves to amplify the audio signal in the lecture~hall.

\textbf{Bidirectional communication.} 
Efficient synchronous bidirectional communication is necessary to engage online students in the lesson. The~use of an airlink connecting the laptop of the professor with the sound system of the lecture hall allows a direct participation and communication not only with the professor but also with on-site students. For~small rooms lacking  a sound system a speakerphone is employed. As~a result online students are actively involved in the lesson, being able to ask questions and participate in open~discussions.

\textbf{Recorded lectures.}
As a backup solution and to mitigate internet connection problems and difficulties related to different time zones, most lectures were also recorded using ZOOM. They were then made available to the students using the online teaching platform Moodle plus Kaltura video cloud service. As~is reflected in the results of the questionnaires described in Section~\ref{sec.evaluation}, the~possibility of recording lectures may lead to an added value to the learning experience by also allowing  the students attending in a synchronous way to review particular parts of the lectures when autonomously studying the~subject.

\subsection{Lecture~Setup}
The setup described above of the on-site equipment was designed aiming at minimizing the complexity of the connections from the professor's point of view. When arriving in the lecture hall the laptop of the professor needs to be physically connected only with a single input/output USB/A cable (coming from the USB hub) and one HDMI cable (towards the HDMI splitter). Then, the~sound system and the ear speaker are paired via bluetooth. Once all connections are established the ZOOM meeting is opened, selecting the correct input and output sound devices and the \textit{share screen} mode is activated. At~this point the professor should select the methodology he/she would like to employ for the lecture:

\begin{description} 
	\item[I.] \textbf{Blackboard}. In~this case the auxiliary video capturing software is opened (VLC, OBS, Quicktime), the~camcorder is selected as the input source and the regular video stream of ZOOM must be turned~off.
	
	\item[II.] \textbf{Graphics tablet}. OpenBoard or the preferred interactive whiteboard software is employed. The~students in the lecture hall will be able to see the screen of the tablet on the projector~screen.
	
	\item[III.] \textbf{Third party software}. While sharing the screen any third party software, like classical mathematical programs such as Matlab, Maple, R Studio, Mathematica, or~PDF and Powerpoint files can be shown. 
\end{description}

Since all three possibilities rely on \textit{screen sharing} via ZOOM and all devices are simultaneously connected, changing between the different options becomes straightforward.

\subsection{Training and Technical~Support}

Despite the fact that the lecture setup described above aimed at being as simple as possible from the professor's point of view,
the blended teaching concept may happen to be a complete novelty for the teaching and technical staff and for the~students.

Thus, already \textit{before the start} of the lectures, we organized a series of special \textit{training sessions}, in~which the professors were trained to use the necessary hardware and software 
in order to get ready for blended teaching before the beginning of the semester. 
During the first part of these sessions, the~new equipment was presented via demonstrations of how each of the three possible lecture modalities (namely the use of blackboard, graphics tablet or third party software)  practically functions.
Moreover, in~order to allow the students to get a preview of how the blended teaching will be in the first semester of 2020/2021,  they were also invited to join some of these training sessions remotely via ZOOM
and to ask any questions they had.	
The presence of hundreds of online students during the demonstrations allowed to show and convince the students and professors of the effectiveness of the adopted concept, including the bidirectional communication.
In the second part of the training sessions, we gave the possibility to each professor (in turn) to personally test the new equipment in order to really familiarize themselves with it, verify the compatibility of their personal laptop with the standard hardware of the lecture halls and be aware of the possibly needed adapters. This was also an occasion for training the technical staff and a special group of students that was hired for adjoint help during the~semester.

Indeed, to~allow blended lectures to start and proceed smoothly,  
a \textit{continuous technical support} was organized \textit{throughout the entire semester}, 
operated by the technical staff of DICAM and by this group of specially trained students.	
During the first week, we offered a \textit{top--down} supporting service, 
with a person of our staff present at the beginning of each lecture in each room ready to facilitate the lecture setup. 
During the rest of the semester, due to the high level of autonomy reached by each professor, the~offered service was changed to the \textit{bottom--up} or \textit{on-demand} type, i.e.,~only in case of necessity the professor could request the intervention of the technical staff always present in the~building.

\section{Evaluation of the Concept and of Its Practical~Implementation} 
\label{sec.evaluation} 

Here we present quantitative results of the evaluation of the blended teaching concept illustrated in the previous sections. {The evaluation of the concept was carried out in three different manners: (1) a technical quality assessment was performed via a direct comparison of the streamed blackboard content using either the normal video streaming feature of ZOOM or~the high-quality low-bandwidth screen sharing strategy proposed in this paper; (2) the authors have carried out quantitative measurements of the required internet bandwidth of the proposed approach; (3) the level of satisfaction with the blended teaching concept proposed in this paper was evaluated separately by the professors and by the students via the specific ex post online questionnaires (PP1) and (PS1).}   
The concept presented in this paper was applied to the entire DICAM department, covering more than 50 courses and 18 lecture halls. 
% T1A, T1B, T2, T3, T4, 1R, EALAB, 2M, 2N, 2A, 2F, 4A 4B, BIB = 14
% 4 for tutor

\subsection{Quantitative Verification of the Video Quality of the Low-Bandwidth High-Quality Screen Sharing Approach Proposed in This Paper Compared to Standard Video~Streaming}
\begin{figure}[]
	%	\begin{center} 
	\includegraphics[width=0.95\textwidth]{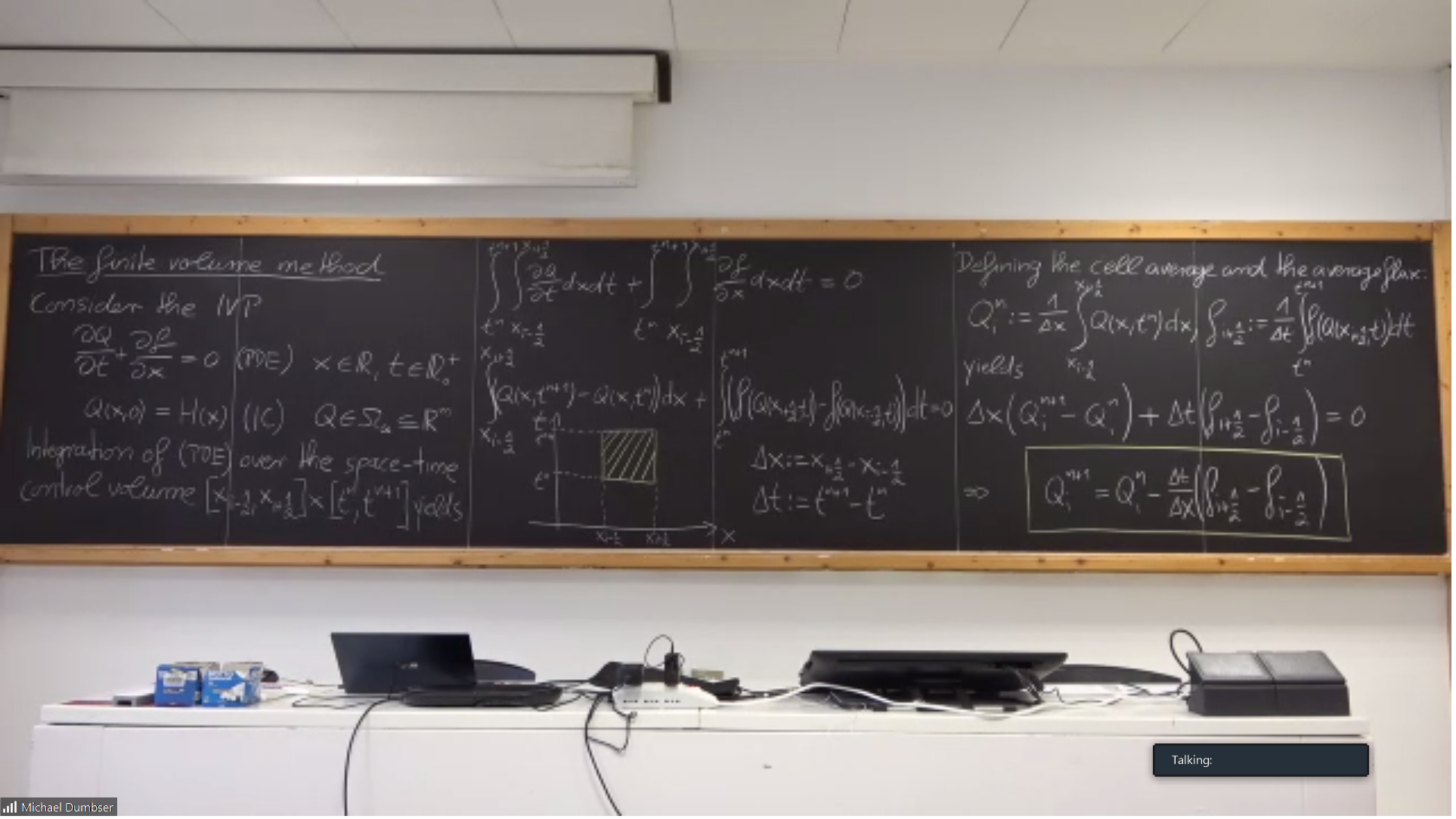} 		
	%	\end{center} 
	\caption{Screenshot of the video content received by a ZOOM client when using the standard ZOOM video feature, which automatically reduces video quality to 360~p due to pandemic internet bandwidth restrictions for meetings with more than two participants. The~blackboard content is only barely~readable. } 	 % AUTHORS to MDPI: this figure and the folowing MUST be shown together in one unique page	
	\label{fig.zoomvideo} 
\end{figure}

\begin{figure}[]	
	%	\begin{center} 
	\includegraphics[width=0.95\textwidth]{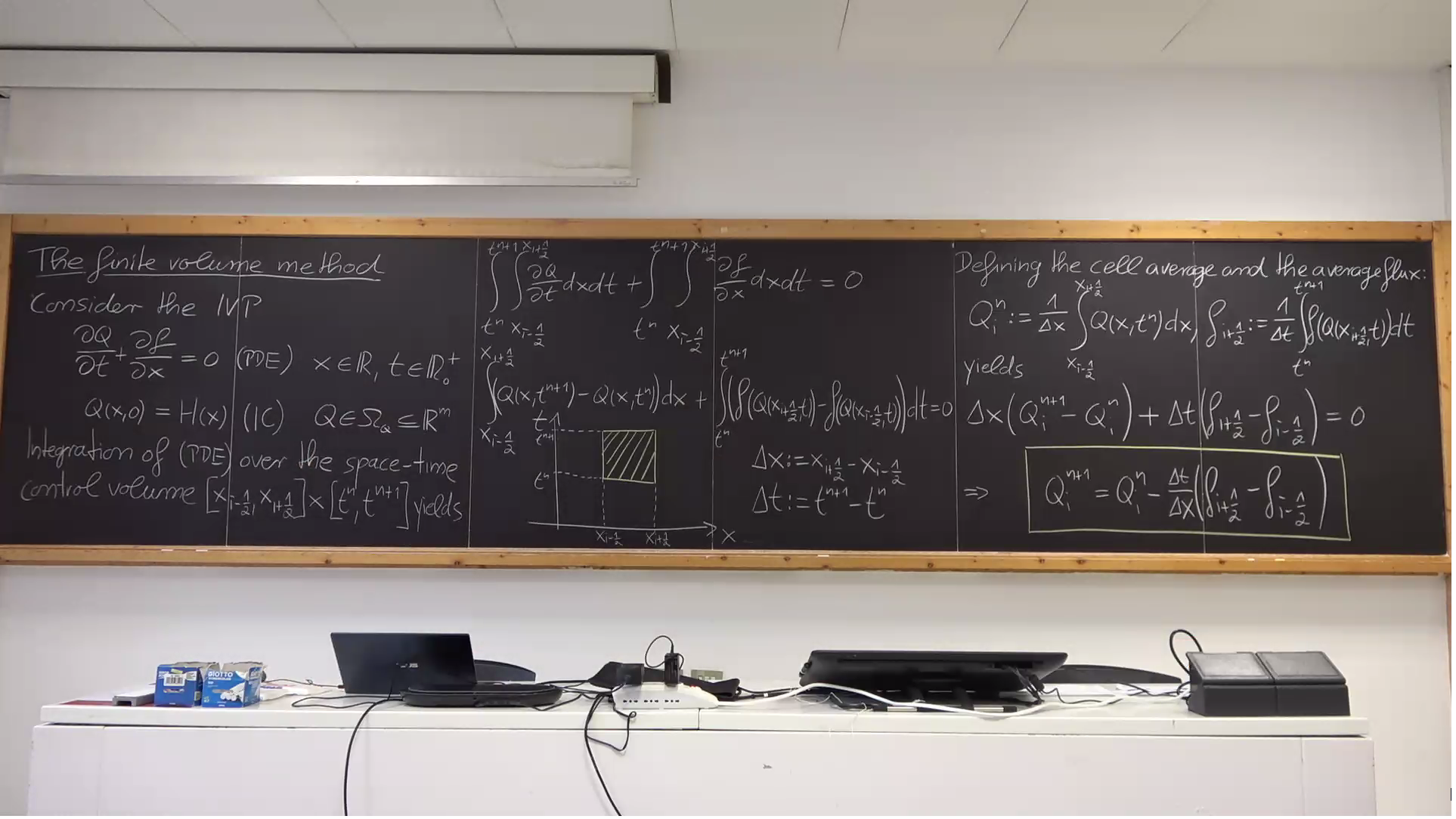} 		
	%	\end{center} 
	\caption{Screenshot of the video content received by a ZOOM client when using the high-quality low-bandwidth screensharing methodology proposed in this paper at 1080~p. The~blackboard content is very clearly readable, even including small indices and~superscripts.} % AUTHORS to MDPI: this figure and the previous one MUST be shown together in one unique page
	\label{fig.screensharing} 
\end{figure}

\begin{figure}[]
	%	\begin{center} 
	\includegraphics[width=0.95\textwidth]{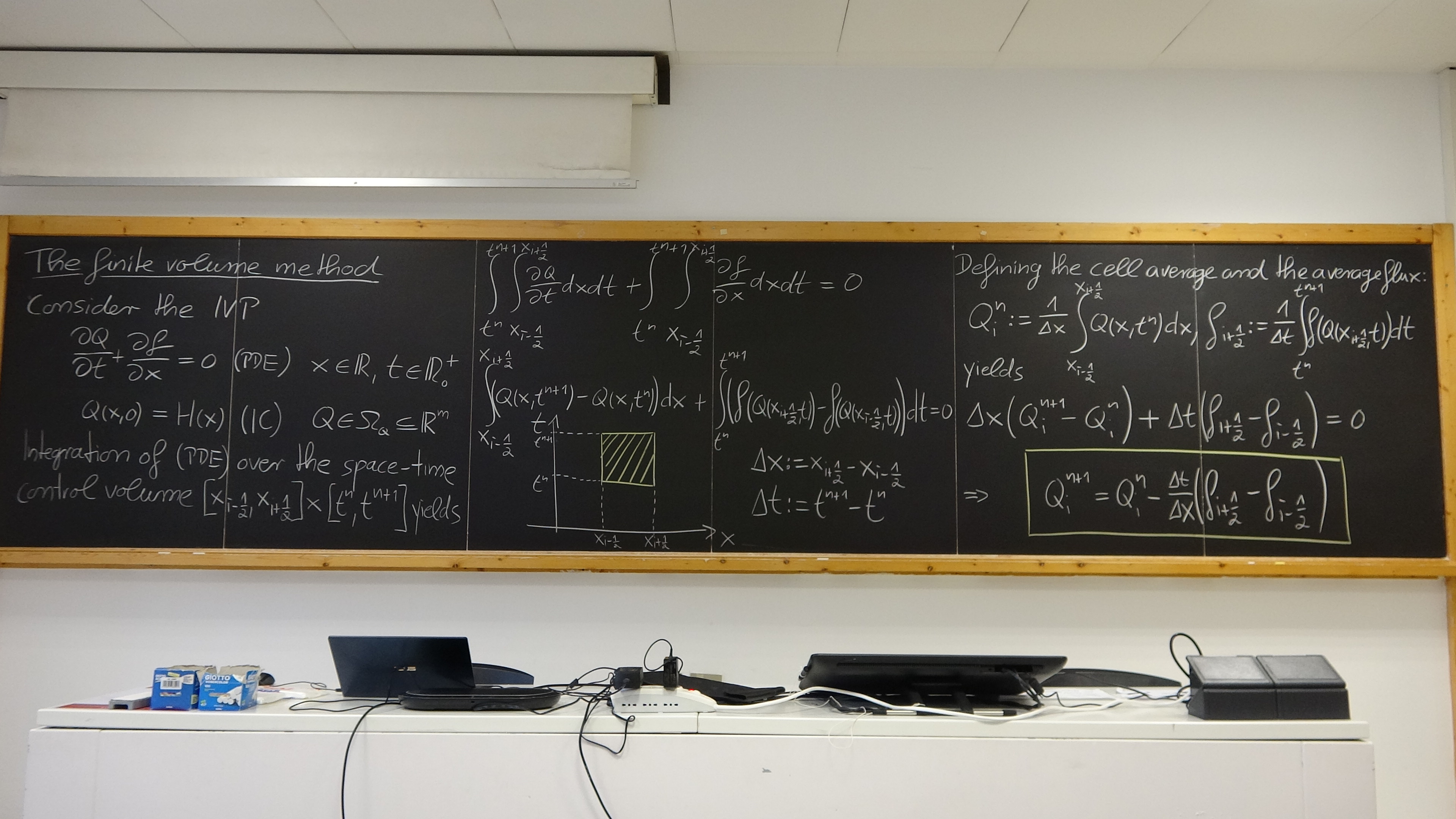} 		
	%	\end{center} 
	\caption{Photograph of the blackboard content taken directly in the lecture hall at 2160p, which serves as a reference  for the previous Figures~\ref{fig.zoomvideo} and ~\ref{fig.screensharing}.} 	 % AUTHORS to MDPI: this figure and the folowing MUST be shown together in one unique page	
	\label{fig.photo}
\end{figure}

In this section we provide a \textit{direct comparison} of the video quality perceived by the students online when using the concept suggested in this paper (full HD screensharing with third party video capturing software at 1080~p) or the standard ZOOM video feature, which automatically reduces the video quality for meetings with more than two participants to 360~p, due to pandemic internet bandwidth restrictions. The~shown blackboard content deliberately contains complex mathematical formulas with small superscripts and indices. In~Figure~\ref{fig.photo} we show a photo of the original blackboard content directly taken in the lecture hall at 2160~p, which serves as a reference. The~screenshots of the transmitted video signals were taken on an Asus Zenbook with full HD screen and ZOOM client 5.3.1 on Windows. For~the video streaming via screen sharing we used ZOOM in combination with OBS Studio 25.0.8. The~video signal was captured with a Sony AX-53 camcorder and an Elgato 4k camlink. In~Figure~\ref{fig.zoomvideo} we show a screenshot of the transmitted video signal using the standard ZOOM video feature, which was automatically limited to only 360~p resolution due to pandemic internet bandwidth restrictions. In~Figure~\ref{fig.screensharing}, instead, we show the transmitted video signal that can be obtained via the proposed low-bandwidth screensharing approach presented in this paper. 
From Figures~\ref{fig.zoomvideo} and~\ref{fig.screensharing} we can clearly conclude that the concept of full HD video streaming via screen sharing proposed in this paper provides a significantly improved video quality compared to the standard ZOOM video feature. In~Figure~\ref{fig.screensharing} the complex mathematical blackboard content is clearly readable, including small subscripts and superscripts, while in Figure~\ref{fig.zoomvideo} it is not. As~already mentioned before, Figure~\ref{fig.photo} shows a photograph directly taken in the lecture hall, which serves as a reference to assess the quality of Figures~\ref{fig.zoomvideo} and~\ref{fig.screensharing}.

\clearpage

	{
		\subsection{Quantitative Internet Bandwidth~Measurements} 	 %Are bold and italic necessary?
		A series of experiments with three different types of internet connections in different locations {{(the three different internet connection types used for the average bandwidth measurements were the following: (C1) direct 1Gbit Ethernet cable connection inside the university building; (C2) ADSL internet connection to the lecture held in Trento from the neighboring region of Lombardy in Northern Italy; (C3) 4G mobile internet connection)}} was carried out. The~obtained average framerates and the necessary internet bandwidth were directly obtained from the \texttt{{statistics feature}} of the various ZOOM clients connected to the lecture. The~measurements showed that with our strategy of video streaming via screen sharing the \textit{{average}} framerate of the video signal of a traditional blackboard lecture transmitted online via ZOOM was still about \textbf{{24 fps}} at permanent full HD resolution, while the total bandwidth required during download by the online students was on average only \textbf{{700 kbps}}, and~hence fully adequate even for low-bandwidth ADSL and average 4G connections. The~reported numbers are arithmetic averages over the three different internet connections C1, C2 and C3 that were tested. These results confirm that the adopted strategy is adequate for the high-quality online streaming of mathematics lectures for most students with the typical internet connections that are nowadays available.       
	}

	\subsection{Quantitative Evaluation of the Blended Teaching Concept via Ex~Post Online Questionnaires}   
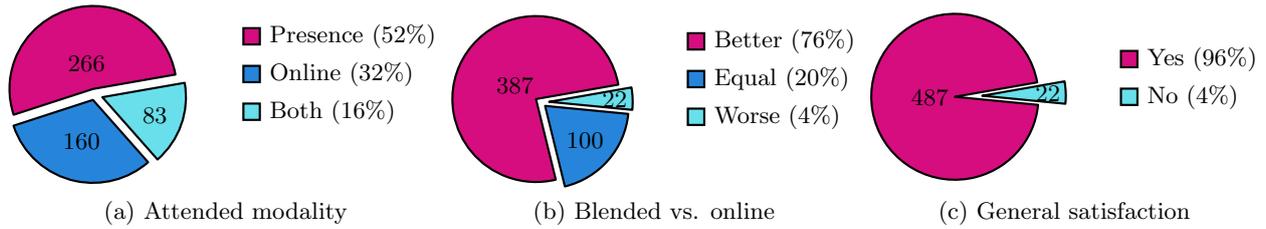
\begin{figure}[!t] \small
	\subfloat[Attended modality]{\begin{tikzpicture}  \footnotesize
		\pie [rotate = 10, explode = {0.075,0.075,0.11}, radius = 1.1, color = {color1,color2,color3},text = legend, sum=509] 
		{ 
			266/Presence  (52\%), 
			160/Online (32\%),
			83/Both (16\%)				
		}
		\end{tikzpicture}}
	\subfloat[Blended vs. online]{\begin{tikzpicture}  \footnotesize
		\pie [rotate = 10, explode = {0.075,0.075,0.11}, radius = 1.1, color = {color1,color2,color3},text = legend, sum=509] 
		%sum=100, sum=100,after number = { percent }, text = inside,cloud, text=inside, polar, explode=0.1, rotate = -65,,text = legend
		{
			387/Better  (76\%), 
			100/Equal (20\%),
			22/Worse (4\%) 
		}
		\end{tikzpicture}}
	\subfloat[General satisfaction]{\begin{tikzpicture}  \footnotesize
		% general satisfaction
		\pie [rotate = 10, explode = 0.18, radius = 1.1, color = {color1,color3},text = legend, sum=509] 
		{
			487/Yes (96\%),  
			22/No (4\%)
		}
		\end{tikzpicture}}
	\caption{Results of the ex post student questionnaire (PS1) concerning the opinion of students  on blended teaching versus pure online teaching (\textbf{b}) and general satisfaction for our blended teaching implementation (\textbf{c}), from~students attending in presence and/or online according to (\textbf{a}). Total answers received:  {$n=509$}. The absolute values of received answers per each option are reported inside the graphs, and~approximate percentages are given in the legends.} % Authors: plese keep the three figures in just one row.
	\label{graph.student.generalfeeling}
\end{figure}

\begin{figure}[!t] \small
	\subfloat[Video quality]{\begin{tikzpicture}  \footnotesize
		\def\printonlybig#1{\ifnum#1>10 #1 \fi}
		\pie [rotate = 10, explode = {0.1,0.05,0.16,0.16}, radius = 1.1, color = {color1,color2,color3,color4},text = legend, sum=509,before number=\printonlybig] 
		{
			256/Very good  (50\%), 
			179/Good (35\%),
			66/Fair (13\%) , 
			8/Not acceptable (2\%)
		}
		\end{tikzpicture}}\qquad \qquad 
	\subfloat[Audio quality]{\begin{tikzpicture} \footnotesize
		\def\printonlybig#1{\ifnum#1>10 #1 \fi}
		\pie [ explode = {0.1,0.05,0.16,0.16}, rotate = 10, radius = 1.1, color = {color1,color2,color3,color4},text = legend, sum=509,before number=\printonlybig] 
		{
			269/Very good (53\%),  
			181/Good (35\%), 
			54/Fair (11\%), 
			5/Not acceptable (1\%)
		}
		\end{tikzpicture}}
	\caption{{Results of the ex post student questionnaire (PS1) concerning} video (\textbf{a}) and audio (\textbf{b}) quality of the blended teaching concept presented in this case study. Total answers received: {$n=509$}. {The relevant absolute values of received answers per each option are reported inside the graphs, and~approximate percentages are given in the legends.}}
	\label{graph.videoaudio}
\end{figure}
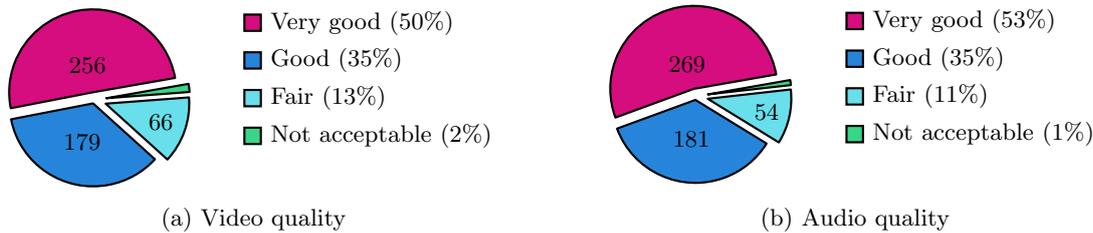

\begin{figure}[!t]
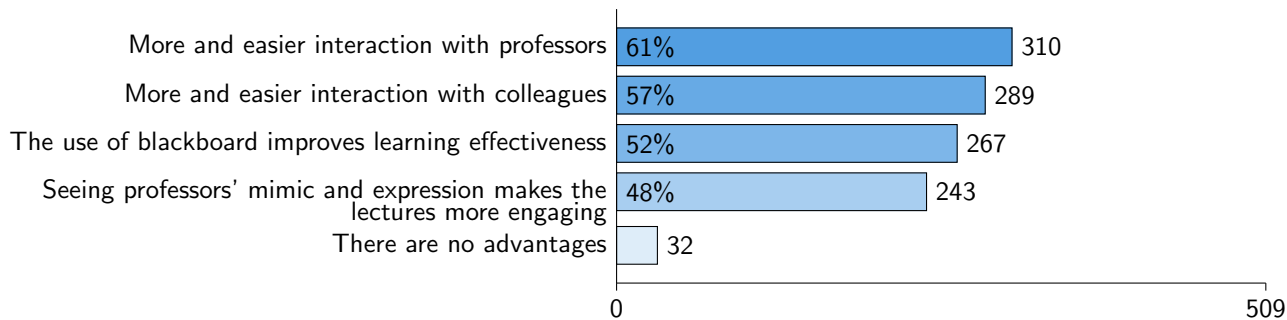
 
	\begin{bchart}[max=509,width=0.52\textwidth] \small
		\bcbar[label=More and easier interaction with professors,text={61\%},color={color2!80}]{310}  %(61\%)
		\bcskip{4pt}
		\bcbar[label=More and easier interaction with colleagues,text={57\%},color=color2!70]{289} %(57\%)
		\bcskip{4pt}
		\bcbar[label=The use of blackboard improves learning effectiveness, text=52\%,color=color2!60]{267}  %(52\%)
		\bcskip{4pt}			
		\bcbar[label=Seeing professors' mimic and expression makes the,text=48\%,color=color2!40]{243} %(48\%)
		\bcskip[label=  lectures more engaging]{2pt}
		\bcskip{4pt}
		\bcbar[label=There are no advantages,color=color2!15]{32}  %(6\%)
	\end{bchart}  
	\caption{New features appreciated in blended teaching vs. pure online teaching {in the ex post student questionnaire (PS1).} Total answers received for the \textit{check box}-type part of the questionnaire: {$n=509$}. 
		{The absolute values of received answers per each option are reported at the end of bars, and~approximate percentages are given inside the bars.}}
	\label{graph.BlendedvsOnline}
\end{figure}

\begin{figure}[!t]
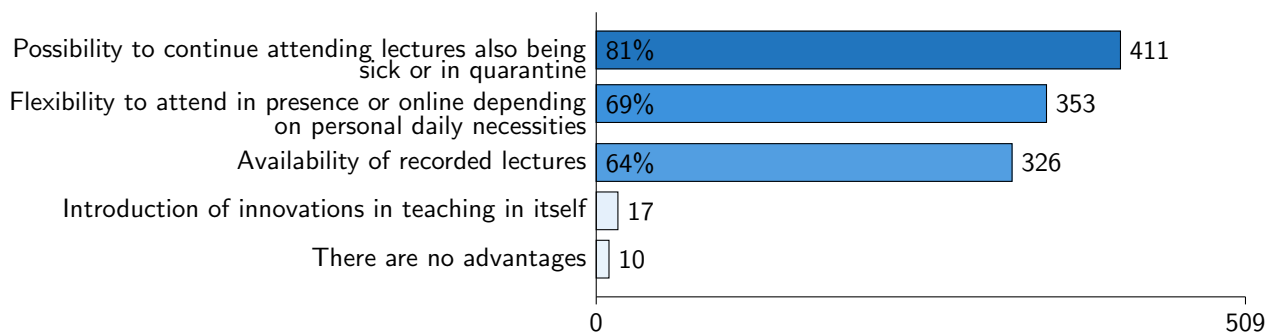

	\begin{bchart}[max=509,width=0.52\textwidth] \small %step =100
		\bcbar[label=\hspace{0.3cm}Possibility to continue attending lectures also being,text=81\%,color=color2!87.5!black]{411}  %(81\%)
		\bcskip[label= sick or in quarantine]{2pt}
		\bcskip[]{4pt}
		\bcbar[label=Flexibility to attend in presence or online depending,text=69\%,color={color2!90}]{353} %(69\%)
		\bcskip[label=on personal daily necessities]{4pt}
		\bcskip[]{4pt}
		\bcbar[label=Availability of recorded lectures,text=64\%,color=color2!80]{326} %(64\%)
		\bcskip{4pt}	
		\bcbar[label=Introduction of innovations in teaching in itself,color=color2!12]{17} % (3\%)
		\bcskip{4pt}
		\bcbar[label=There are no advantages,color=color2!10]{10} % (2\%)	
	\end{bchart} 
	\caption{New features appreciated in blended teaching vs. traditional teaching {in the ex post student questionnaire (PS1).} Total answers received for the \textit{check~box}-type part of the questionnaire: {$n=509$}. 
		{The absolute values of received answers for each option are reported at the end of bars, and~approximate percentages are given inside the bars.}}
	\label{graph.BlendedvsTraditional}
\end{figure}

	For the quantitative evaluation of the level of appreciation of the blended teaching concept proposed in this paper, we have considered all of the courses involving teaching of mathematics and its application to theoretical mechanics; we aimed for the evaluation to be representative of each year of study, %please confirm rewording
	and thus examined the courses of \textit{Analysis 1} (BSc, first year), \textit{Analysis 2} (BSc, second year), \textit{Theoretical solid mechanics} (BSc, third year), \textit{Theoretical fluid mechanics} (BSc, third year) and \textit{Numerical analysis for PDE} (MSc, first year). For~the ex post questionnaire sent to the students (PS1), we received a total of $n = 509$ answers. Meanwhile the opinions of the professors were collected via the ex post questionnaire (PP1) that obtained $n=6$ answers from the professors leading the didactic team of each course. The~questionnaires were compiled after one month of blended lectures. The~sample size of $n=6$ concerning the answers received for the ex post questionnaire (PP1) was small since the objective was to focus on the aspects of blended teaching in mathematics and its application to theoretical mechanics, where almost all lectures are traditionally held at the blackboard.   
	The questions and answers to the ex post student questionnaire (PS1) together with the related quantitative results can be found in Table~\ref{tab.ps1.survey} and are also graphically presented in Figures~\ref{graph.student.generalfeeling}--\ref{graph.BlendedvsTraditional}. 
	
\begin{table}[]
	\caption{{\textit{Ex post} survey (PS1) on the blended teaching concept of DICAM made among students. Total number of answers received $n=509$. } }
	\label{tab.ps1.survey} 
{ \small
		\begin{tabular}{l|l} 
			\hline 
			\textbf{Questions} & \textbf{Answers with absolute numbers} \\ 
			\hline 
			1. How do you typically attend the lectures? &   \phantom{1 -} \textbf{in presence (266)} \\
			&  \phantom{1 -} online (160) \\
			&  \phantom{1 -} both (online and in presence) (83) \\ 
			\hline 
			2. How do you evaluate the quality of the blended teaching  &  \phantom{1 -} \textbf{better (387)} \\ 
			\hspace{2.3mm} compared to pure online teaching?            &  \phantom{1 -} equal (100) \\ 
			&  \phantom{1 -} worse (22) \\ 
			\hline 
			3. Are you overall satisfied with the blended teaching     &   \phantom{1 -} \textbf{yes (487)} \\ 
			\hspace{2.3mm} concept offered at DICAM?                   &   \phantom{1 -} no (22) \\ 
			\hline 
			4. How do you evaluate the video quality of the blended teaching     &    \textbf{4 - very good (256)} \\ 
			\hspace{2.3mm} 	concept at DICAM? &  3 - good (179) \\ 
			&  2 - fair (66) \\ 
			&  1 - not acceptable (8) \\ 
			\hline 
			5. How do you evaluate the audio quality of the blended teaching     &    \textbf{4 - very good (269)} \\ 
			\hspace{2.3mm} concept at DICAM?                                                    &  3 - good (181) \\ 
			&  2 - fair (54) \\ 
			&  1 - not acceptable (5) \\ 
			\hline 
		\end{tabular}
	}
	{\small
		\begin{tabular}{l|l} 
			6. Which are the main advantages of             &   \textbf{more and easier interaction with professors (310) } \\ 
			\hspace{2.3mm} blended teaching compared to                                         &  more and easier interaction with colleagues (289) \\ 
			\hspace{2.3mm} pure online teaching? &  use of blackboard improves learning efficiency (267) \\ 
			\hspace{2.3mm}  (multiple answers allowed) &  seeing the professor makes lectures more engaging (243)  \\    
			&  there are no advantages (32) \\ 	
			\hline 
			7. Which are the main advantages of             &   \textbf{possibility to attend lectures also being sick/in quarantine (411) } \\ 
			\hspace{2.3mm} blended teaching compared to                                         &  flexibility to attend in presence or online acc. to personal needs (353) \\ 
			\hspace{2.3mm} traditional teaching?  &  availability of recorded lectures (326) \\ 
			\hspace{2.3mm} (multiple answers allowed) &  introduction of innovations in teaching in itself (17)  \\    
			&  there are no advantages (10) \\ 	
			\hline 			
		\end{tabular} 
	} 
\end{table}

The first set of questions of (PS1) was about the way in which students participate in the lectures, how they evaluate the quality of blended teaching versus pure online teaching and whether the students were globally satisfied with the concept presented in this paper or~not; see Table~\ref{tab.ps1.survey} and Figure~\ref{graph.student.generalfeeling}.   
More than half of the total answers comes from students attending the courses mostly in presence, 
a third of them attended mostly online and a part of them chose the preferred modality on a daily base depending on personal necessities. Almost the totality of the received answers shows a \textit{general satisfaction} with the implemented blended teaching concept and a \textit{clear preference} for blended teaching with respect to pure online teaching; see the detailed results reported in Figure~\ref{graph.student.generalfeeling}.

The second set of questions was specifically designed to assess the audio and video quality perceived by the students; see  Table~\ref{tab.ps1.survey}  and Figure~\ref{graph.videoaudio}.     
Indeed, the~overall video and audio quality was very positively judged by around 85\% of the students, who attributed a score of $\geq 3$, i.e.,~\textit{good} or \textit{very good}. 
{In order to provide a clearer analysis of these two fundamental questions in the (PS1) questionnaire, 
	we have assigned a quantitative score from $1$ to $4$ to each qualitative judgment received from the students, 
	i.e.,~$1$ to \textit{Not acceptable}, $2$ to \textit{Fair}, $3$ to \textit{Good}, and~$4$ to \textit{Very good}; see
	also Table~\ref{tab.ps1.survey}.  
	From the collected data we can state that the video quality was evaluated with a mean score of $3.34$ out of $4$, with~a standard deviation of $\sigma=0.76$. The~audio quality obtained a mean score of $3.40$ out of $4$, with~a standard deviation of $\sigma=0.72$. See also Table~\ref{tab.ps1.survey} and Figure~\ref{graph.videoaudio} for more details about the exact numbers concerning the collected answers.}

Moreover, the~key assumptions guiding the conception of our blended teaching approach were confirmed by the opinions of DICAM students: see the last two questions in Table~\ref{tab.ps1.survey} and Figures~\ref{graph.BlendedvsOnline} and  \ref{graph.BlendedvsTraditional}, where we summarize the results of \textit{check box}-type questionnaires asking which new features introduced with blended teaching were important for the students. 
With respect to pure online teaching, students appeared to appreciate the major and easier interaction with professors and colleagues and the increased effectiveness of lectures thanks to the use of the blackboard and the possibility to easily see professors' body language and expression. 
With respect to standard traditional teaching, they appreciated the new high flexibility of switching between in-presence and online lectures, not only for COVID-19-related emergency reasons, but~also for normal everyday life organization. 
The availability of recorded lectures was greatly appreciated also in the free comment space of the questionnaire.


 {A quantitative analysis of the opinion of the professors from the area of mathematics and its application to theoretical mechanics was carried out with the aid of the ex post survey (PP1), the~questions and received answers of which are reported in Table~\ref{tab.pp1.survey}. From~the received answers, it becomes clear that professors of mathematics and theoretical mechanics} show a clear preference for the blended teaching concept with respect to  pure online teaching; in particular, they highly appreciated the possibility to continue using traditional blackboards, showing experiments and having a true interaction with the students physically present in the lecture hall.
Some of them were also willing to maintain several of the concepts  introduced here even when traditional teaching will be again possible, considering indeed that the possibility of recording, the~increased flexibility and the new available technological instruments can provide an added value to the lectures and widen the reached~target.

\begin{table}[!ht]
	\caption{\textcolor{black}{\textit{Ex post} survey (PP1) of DICAM made among its professors  teaching mathematics and theoretical mechanics for civil and environmental engineering students. Total number of answers received $n=6$. } }
	\label{tab.pp1.survey} 
	\textcolor{black}{ \small 
		\begin{tabular}{l|l} 
			\hline 
			\textbf{Questions} & \textbf{Answers with absolute numbers} \\ 
			\hline 
			1. Which instruments do you use for teaching? &  \phantom{1 -} \textbf{blackboard (5)} \\
			&  \phantom{1 -} graphics tablet (1) \\
			&  \phantom{1 -} slides (0) \\ 
			\hline 
			2. How do you evaluate the quality of the technical equipment  &   \textbf{4 - excellent (4)} \\ 
			\hspace{2.3mm} made available by DICAM for blended teaching?                  &  3 - very good (2) \\ 
			&  2 - good (0) \\ 
			&  1 - fair (0) \\ 
			&  0 - insufficient (0) \\ 
			\hline 
			3. How do you evaluate the quality of the training courses     &   \textbf{4 - excellent (4)} \\ 
			\hspace{2.3mm} offered by DICAM for blended teaching?          &  3 - very good (2) \\ 
			&  2 - good (0) \\ 
			&  1 - fair (0) \\ 
			&  0 - insufficient (0) \\ 
			\hline 
			4. How important is technical support for blended teaching?    &   \textbf{3 - very important (3)} \\ 
			&  2 - important (2) \\ 
			&  1 - indifferent (1) \\ 
			&  0 - not needed (0) \\ 
			\hline 
			5. How do you evaluate the quality of the technical support    &   \textbf{4 - excellent (5)} \\ 
			\hspace{2.3mm} provided by DICAM for blended teaching?         &  3 - very good (1) \\ 
			&  2 - good (0) \\ 
			&  1 - fair (0) \\ 
			&  0 - insufficient (0) \\ 
			\hline 
			6. How do you evaluate the quality of bidirectional communication    &   \textbf{4 - excellent (1)} \\ 
			\hspace{2.3mm} with the students online?                       &  3 - very good (3) \\ 
			&  2 - good (0) \\ 
			&  1 - fair (2) \\ 
			&  0 - insufficient (0) \\ 
			\hline 
			7. Compared to \textit{pure online teaching}, how do you evaluate the   &   \textbf{3 -much better (4)} \\ 
			\hspace{2.3mm} overall quality of \textit{onsite blended teaching}?    &  2 - better (2) \\ 
			&  1 - the same (0) \\ 
			&  0 - worse (0) \\ 
			\hline 
		\end{tabular}
	}
	\textcolor{black}{ \small
		\begin{tabular}{l|l} 
			8. Which are the main advantages of        &   \textbf{the possibility to use traditional teaching instruments (4) } \\ 
			\hspace{3.0mm} blended teaching compared      &  the possibility to see the students' reactions (3) \\ 
			\hspace{3.0mm} to pure online teaching?  &  improves the direct interaction with the students (2) \\ 
			\hspace{3.0mm} (multiple answers allowed) &  \textbf{improves the overall quality of teaching (4)} \\    	
			\hline 
		\end{tabular} 
	} 
\end{table}

\subsection{Moving Beyond: Innovative Didactic Concepts for the Numerical Analysis~Course}

Given the available new technological equipment and in order both (i) to mitigate the adjoint difficulties students can perceive due to the pandemic crisis, but~also (ii) to provide additional and high-quality learning opportunities, 
a specific teaching concept was designed and tested for the course \textit{Numerical analysis for PDE} (MSc, first year) held in 2020/2021 by the first two authors of this paper. For~the contents of the course we refer to~\cite{ValliQuarteroni,Toro}.
In addition to the blended theoretical lectures at the blackboard, according to the technical concept outlined in this paper, the~course offered: (i) weekly synchronous online computer exercises with from-scratch implementation of the studied numerical methods; (ii) weekly exercise sheets with reference solutions provided under the form of a \textit{podcast/screencast}, i.e.,~audio and video recordings of the commented solutions of the exercise sheets; (iii) two projects concerning the numerical solution of simplified real life applications, to~be developed in small groups and with students documenting their findings in a written report and presenting and discussing the methods and results in class, following the concepts of a \textit{continuous assessment} of the students~\cite{sanz2019students} and the \textit{flipped classroom principle}~\cite{CFW17,Cong2020}; (iv) mixed use of Italian and English language to also transmit  specific technical vocabulary and the capacity to convey subject-related concepts at an international~level. 

After five weeks of lectures, corresponding to almost the 50\% of the entire course, the~25 students answering the questionnaire were overall satisfied with this new offer, all of them found  the weakly exercise sheets and the commented solutions in the form of podcast to be useful, and~22 students appreciated the mixed language approach, with 18 of them thinking this results in an added value for the course. The~project to be developed by the students was elaborated with mixed groups in which some of the students were online while others were physically present in the lecture hall. To~this end a good \textit{bidirectional communication} was crucial. The~technical equipment  allowed the students to interact with each others and with the professor in a fluent and synchronous~way.

\section{Conclusions}
\label{sec.conclusions}

In this paper we have presented all the technical details needed to realize an economically affordable and simple but efficient technical concept for the realization of blended teaching of mathematics for engineering students during the COVID-19 pandemic. {The case study reported in this paper was carried out at the level of an entire engineering department at the University of Trento, Italy, with~$n=1011$ students and $n=68$ professors participating. The~presented case study was carefully designed and guided by ex ante questionnaires sent to students and professors before the beginning of the semester. } The technical concept that was realized is based on  three key assumptions: (i) that traditional blackboard lectures, including the gestures and facial expressions of the professor, are still a very efficient and highly appreciated means of teaching mathematics at the university; (ii) that a lecture in mathematics is a creative process that requires the physical presence of an audience so that the lecture can be best adjusted to the needs of the students; and finally (iii) that undergraduate students need a minimum level of direct personal interactions and the possibility to get to know each other, especially in the first years. 
All three basic assumptions, as~well as the perception of the audio and video quality, have been quantitatively verified with the aid of systematic ex post online surveys that were sent to all students. The~overall level of satisfaction with the blended teaching concept presented in this paper was about 96\%. A~series of direct measurements revealed that the total internet bandwidth required to download  full HD video streaming via screen sharing with ZOOM was rather low (only about 700 kbps) and thus makes the methodology applicable also during pandemic internet bandwidth restrictions and makes blended teaching also available for those students who have a fairly slow internet connection. If~necessary and if possible, this concept could be adopted %confirm "adopted".
also in subsequent semesters {during the COVID-19 pandemic. While the present paper mainly focuses on the problem of blended teaching of \textit{mathematics lectures}, which usually exhibit  complex blackboard content with complicated formulas containing small subscripts and superscripts, the~same concept can also be easily transferred to other STEM subjects with complex blackboard contents, like physics, chemistry or biology.
	Given the low internet bandwidth requirements and the affordable cost of the concept presented in this paper, the~authors believe that an implementation in other cultural contexts is also possible, but~this extension is clearly beyond the scope of this paper, which is limited to one specific case study.}

%%%%%%%%%%%%%%%%%%%%%%%%%%%%%%%%%%%%%%%%%%
\vspace{6pt} 

%%%%%%%%%%%%%%%%%%%%%%%%%%%%%%%%%%%%%%%%%%
\paragraph{Author Contributions.}{Conceptualization, methodology and implementation: S.B., M.D., E.G. Writing original draft and editing: S.B., M.D., E.G.; Funding Acquisition: M.D.. }

%%%%%%%%%%%%%%%%%%%%%%%%%%%%%%%%%%%%%%%%%%
\paragraph{Funding.} {The technical equipment needed for the blended teaching concept detailed in this paper was funded by the Department of Civil and Environmental Engineering (DICAM) of the University of Trento and by the University of Trento, Italy. }

%%%%%%%%%%%%%%%%%%%%%%%%%%%%%%%%%%%%%%%%%%
\paragraph{Acknowledgments.} {
	The authors would like to thank the following professors, students and technical staff of University of Trento for the technical and logistic support that was necessary for the realization of the blended teaching concept proposed in this paper. 
	Professors: Prof. Dr. A. Vitti; technical staff of UniTN: S. Bernardini, C. Casagranda, I. Cristofolini, M. Filippi, M. Paoletto and F. Romagnoli; students: 
	G. Bagattini, M. Bassetti, M. Bressanin, C. Cassar\`a, M. Dalpiaz, M. Facchin, F. Leali, 
	A. Lonardi, R. Panziera, C. Paoli, G. Scarcella, L. Trentini and D. Vallenari.   
	
	The authors also would like to thank all the professors in charge of the courses in analysis and theoretical solid and fluid mechanics for their kind support and patience: Prof. F. Bagagiolo, Prof. D. Bigoni, Prof. L. Fraccarrollo, Prof. N. Pugno, Prof. G. Rosatti and Prof. A. Valli. 
	
	Special thanks also to Prof. Dr. Oreste S. Bursi, Director of the Department of Civil, Environmental and Mechanical Engineering, for granting all the necessary support and freedom to realize the concept detailed in this paper. 
	
	\textcolor{black}{Last, but not least, the authors would also like to kindly thank the two anonymous referees for their very useful and constructive comments, which helped to improve the quality of this paper substantially. }
}

%%%%%%%%%%%%%%%%%%%%%%%%%%%%%%%%%%%%%%%%%%
\paragraph{Conflicts of interest.}{The authors declare no conflict of interest.} 

\paragraph{Authors:}
\begin{description}
	\item[Saray Busto] obtained her PhD degree in mathematics \textit{summa cum laude} in 2018 at the University of Santiago de Compostela, Spain. She is working in the field of numerical analysis for hyperbolic partial differential equations. 
	\item[Michael Dumbser] obtained his PhD degree in aerospace engineering \textit{summa cum laude} in 2005 at the University of Stuttgart, Germany. He is working in the field of numerical analysis for hyperbolic partial differential equations. Since 2018 he is Dean of Studies of the Department of Civil, Environmental and Mechanical Engineering (DICAM) of the University of Trento, Italy. Together with C.D. Munz and S. Roller he was winner of the teaching award 2003 of the Ministry of University and Research of the State of Baden-W\"urttemberg, Germany, for the lecture \textit{Numerical gasdynamics} held at the University of Stuttgart, Germany, in 2003.
	\item[Elena Gaburro] obtained her PhD degree in mathematics \textit{summa cum laude} in 2018 at the University of Trento, Italy. She is working in the field of numerical analysis for hyperbolic partial differential equations.  	
\end{description}

\bibliographystyle{plain}

\end{document}